\documentclass[preprint,11pt,times]{elsarticle}

\usepackage{setspace}
\usepackage{geometry}
\usepackage{url}
\usepackage{graphicx,graphics}
\usepackage{subfigure}
\usepackage{epsfig}
\usepackage{psfrag}
\usepackage{hyperref}
\usepackage{amsmath,amssymb,amsfonts,amsthm}
\usepackage{mathrsfs}
\usepackage{color}
\usepackage{float}
\usepackage{natbib}
\usepackage{tabularx}
\usepackage{multirow}

\usepackage{pgfplots}
\usetikzlibrary{plotmarks}

\theoremstyle{remark}

\newcommand{\Eref}[1]{Equation (\ref{#1})}
\newcommand{\fref}[1]{Figure (\ref{#1})}
\newcommand{\Erefs}[1]{Equations (\ref{#1})}
\newcommand{\frefs}[1]{Figures~\ref{#1}}

\newcommand{\rmd}{\mathrm{d}}
\newcommand{\bveps}{\boldsymbol{\varepsilon}}
\newcommand{\bvsig}{\boldsymbol{\sigma}}

\newcommand{\cc}{\mathbf{C}}
\newcommand{\dd}{\mathbf{D}}

\newcommand{\KK}{\mathbf{K}}
\newcommand{\kk}{\mathbf{K}}
\newcommand{\bfm}{\mathbf{M}}

\newcommand{\qq}{\mathbf{q}}

\newcommand{\uu}{\mathbf{u}}

\begin{document}

\begin{frontmatter}

\title{Supersonic flutter analysis of flat composite panels by unified formulation}       

\author[unsw]{S Natarajan \corref{cor1}\fnref{fn1}}
\author[india]{Ganapathi Manickam}
\author[portugal,dubai]{AJM Ferreira}
\author[italy,dubai]{E Carrera}

\cortext[cor1]{Corresponding author}

\address[unsw]{School of Civil \& Environmental Engineering, The University of New South Wales, Sydney, Australia}
\address[india]{Head, Stress \& DTA, IES-Aerospace, Mahindra Satyam Computers Services Ltd., Bangalore, India}
\address[portugal]{Faculdade de Engenharia da Universidade do Porto, Porto, Portugal.}
\address[italy]{Department of Aeronautics and Aerospace Engineering, Politecnico di Torino, Italy.}
\address[dubai]{Department of Mathematics, Faculty of Science, King Abdulaziz University, P.O. Box 80203, Jeddah 21589, Saudi Arabia.}
\fntext[fn1]{School of Civil \& Environmental Engineering, The University of New South Wales, Sydney, NSW 2052, Australia. Tel:+61 293855030, Email: s.natarajan@unsw.edu.au; sundararajan.natarajan@gmail.com}

\begin{abstract}
In this paper, the linear flutter characteristics of laminated composite flat panels immersed in a supersonic flow is investigated using field consistent elements within the framework of unified formulation. The influence of the aerodynamic damping on the supersonic flutter characteristics of flat composite panels is also investigated. The aerodynamic force is evaluated using two-dimensional static aerodynamic approximation for high supersonic flow. Numerical results are presented for laminated composites that bring out the influence of the flow angle, the boundary conditions, the plate thickness and the plate aspect ratio on the flutter characteristics.
\end{abstract}

\begin{keyword} 
	Laminated plates \sep unified formulation \sep supersonic flutter analysis \sep field consistent element \sep shear locking \sep sinusoidal shear deformation theory
\end{keyword}

\end{frontmatter}


\section{Introduction}
Engineered materials, such as the laminated composites are extensively used in various engineering disciplines such as in the aerospace engineering, automotive engineering and civil engineering. The laminated composite structures are often made of several orthotropic layers with different materials stacked together to achieve desired properties. Such a construction provides flexibility in tailoring the properties of the structures by varying the stack up sequence or by changing the fiber orientation in the lamina. In practice, the use of these materials in aerospace industries has necessitated to understand the dynamic characteristics of laminated structures, especially when exposed to air flow. The phenomenon called `panel flutter' involves the interaction of panel inertia force, the elastic restoring force, the thermal force and the air stream which passes over one side of the panel. This study is important in aerospace structural design in evaluating the fatigue life and allowable stress of the structural component exposed to supersonic flow. This has attracted researchers~\cite{singhaganapathi2005, wang1997} to study the flutter characteristics of composite panels. Sawyer~\cite{sawyer1977} using classical lamination theory studied the flutter characteristics of laminated plates. Srinivasan and Babu~\cite{srinivasanbabu1987} examined the flutter of laminated quadrilateral plates. Liao and Sun~\cite{liaosun1993} investigated the supersonic flutter behaviour of stiffened composite skew plates using a degenerated shell element. Pidaparti and Chang~\cite{pidapartichang1998} investigated the flutter characteristics of skewed and cracked composite panels. However, all these investigations neglected the effect of damping. Singha and Ganapathi~\cite{singhaganapathi2005} studied the flutter characteristics of skew composite plates by considering damping and thermo-mechanical loads using a shear deformable high precision 4-noded quadrilateral element.

In describing the plate kinematics, several laminate plate theories have been applied for the analysis of composite plates, such as the classical laminated plate theory (CLPT)~\cite{bosereddy1998}, first order shear deformation theory (FSDT)~\cite{whitneypagano1970} and other higher-order shear deformation theories (HSDTs)~\cite{reddy1984,kantmanjunatha1988,pandyakant1988,ganapathipatel2004}. Recently Carrera~\cite{carrera2001} derived a series of axiomatic approaches, coined as Carrera Unified Formulation~\cite{carrerademasi2002} for the general description of two-dimensional formulations for multilayered plates and shells. With this unified formulation, it is possible to implement in a single software a series of hierarchical formulations, thus affording a systematic assessment of different theories ranging from simple equivalent single layer models up to higher order layerwise descriptions. The plate structures are analyzed by employing a numerical technique. One such powerful and versatile technique is the finite element method. The CUF has been used to develop discrete models such as the finite element method~\cite{carrerademasi2002,carrerademasi2002a} and, more recently, meshless methods based upon collocation with radial basis functions~\cite{ferreiraroque2011}. 

\paragraph{Objective} In this paper, a $\mathcal{C}^o$ 4-noded quadrilateral shear flexible element is employed to study the free vibration and flutter characteristics of laminated composite plates immersed in a supersonic flow. The plate kinematics is based on Carrera Unified Formulation (CUF) and a hybrid displacement assumption is used for the in-plane and the transverse displacements. The shear locking is alleviated by employing a field consistent approach. The influence of the flow angle, the plate thickness, the plate aspect ratio, the boundary conditions and damping on the critical aerodynamic pressure is numerically studied. 

\paragraph{Outline} The paper commences with a brief discussion on the unified formulation for plates and the finite element discretization. Section \ref{edesc} describes the element employed for this study. The efficiency of the present formulation, numerical results and parametric studies are presented in Section~\ref{numres}, followed by concluding remarks in the last section.

\section{Carrera Unified Formulation}\label{cuftheory}
\subsection{Basis of CUF}
Let us consider a laminated plate composed of perfectly bonded layers with coordinates $x,y$ along the in-plane directions and $z$ along the thickness direction of the whole plate, while $z_k$ is the thickness of the $k^{\rm th}$ layer. The CUF is a useful tool to implement a large number of two-dimensional models with the description at the layer level as the starting point. By following the axiomatic modelling approach, the displacements $\uu(x,y,z) = ( u(x,y,z), v(x,y,z), w(x,y,z))$ are written according to the general expansion as:
\begin{equation}
\uu(x,y,z) = \sum\limits_{\tau = 0}^N F_\tau(z) \uu_\tau(x,y)
\label{eqn:unifieddisp}
\end{equation}
where $F(z)$ are known functions to model the thickness distribution of the unknowns, $N$ is the order of the expansion assumed for the through-thickness behaviour. By varying the free parameter $N$, a \emph{hierarchical} series of two-dimensional models can be obtained. The strains are related to the displacement field via the geometrical relations:
\begin{eqnarray}
\bveps_{pG} = \left[ \begin{array}{ccc} \varepsilon_{xx} & \varepsilon_{yy} & \gamma_{xy}  \end{array} \right]^{\rm T} = \dd_p \uu \nonumber \\
\bveps_{nG} = \left[ \begin{array}{ccc} \gamma_{xz} & \gamma_{yz} & \varepsilon_{zz} \end{array} \right]^{\rm T} = \left( \dd_{np} + \dd_{nz} \right) \uu
\label{eqn:strainDisp}
\end{eqnarray}
where the subscript $G$ indicate the geometrical equations, $\dd_{p}, \dd_{np}$ and $\dd_{nz}$ are differential operators given by:
\begin{eqnarray}
\dd_p = \left[ \begin{array}{ccc} \partial_x & 0 & 0 \\ 0 & \partial_y & 0 \\ \partial_y & \partial_x & 0 \end{array} \right], \hspace{0.5cm} \dd_{np} = \left[ \begin{array}{ccc} 0 & 0 & \partial_x \\ 0 & 0 & \partial_y \\ 0 & 0 & 0 \end{array} \right], \nonumber \\
\dd_{nz} = \left[ \begin{array}{ccc} \partial_z & 0 & 0 \\ 0 & \partial_z & 0 \\ 0 & 0 & \partial_z \end{array} \right].
\end{eqnarray}
The 3D constitutive equations are given as:
\begin{eqnarray}
\bvsig_{pC} = \cc_{pp} \bveps_{pG} + \cc_{pn} \bveps_{nG} \nonumber \\
\bvsig_{nC} = \cc_{np} \bveps_{pG} + \cc_{nn} \bveps_{nG}
\label{eqn:stressdef}
\end{eqnarray}
with 
\begin{eqnarray}
\cc_{pp} = \left[ \begin{array}{ccc} C_{11} & C_{12} & C_{16} \\ C_{12} & C_{22} & C_{26} \\ C_{16} & C_{26} & C_{66} \end{array} \right] \hspace{0.5cm} \cc_{pn} = \left[ \begin{array}{ccc} 0 & 0 & C_{13} \\ 0 & 0 & C_{23} \\ 0 & 0 & C_{36} \end{array} \right] \nonumber \\
\cc_{np} = \left[ \begin{array}{ccc} 0 & 0 & 0 \\ 0 & 0 & 0 \\ C_{13} & C_{23} & C_{36} \end{array} \right] \hspace{0.5cm} \cc_{nn} = \left[ \begin{array}{ccc} C_{55} & C_{45} & 0 \\ C_{45} & C_{44} & 0 \\ 0 & 0 & C_{33} \end{array} \right]
\end{eqnarray}
where the subscript $C$ indicate the constitutive equations. The \textit{Principle of Virtual Displacements} (PVD) in case of multilayered plate subjected to mechanical loads is written as:
\begin{equation}
\sum\limits_{k=1}^{N_k} \int\limits_{\Omega_k}\int\limits_{A_k} \left\{  (\delta \bveps_{pG}^k)^{\rm T} \bvsig_{pC}^k + (\delta \bveps_{nG}^k)^{\rm T} \bvsig_{nC}^k \right\}~\mathrm{d}\Omega_k~\mathrm{d}z = \sum\limits_{k=1}^{N_k} \int\limits_{\Omega_k} \int\limits_{A_k} \rho^k \delta \uu_s^{k^{\rm T}} \ddot{\uu}^k~\mathrm{d}\Omega_k~\mathrm{d}z + \sum\limits_{k=1}^{N_k} \delta \mathbf{L}_e^k 
\end{equation}
where $\rho^k$ is the mass density of the $k^{\rm th}$ layer, $\Omega_k$, $A_k$ are the integration domain in the $(x,y)$ and the $z$ direction, respectively. Upon substituting the geometric relations (\Eref{eqn:strainDisp}), the constitutive relations (\Eref{eqn:stressdef}) and the unified formulation into the PVD statement, we have:
\begin{equation}
\begin{split}
\int\limits_{\Omega_k}\int\limits_{A_k} \left\{ \left(\dd_p^k F_s \delta \uu_s^k\right)^{\rm T} \left\{ \cc_{pp}^k \dd_p^k F_{\tau} \uu_{\tau}^k + \cc_{pn}^k (\dd_{n\Omega}^k + \dd_{nz}^k) F_{\tau}\uu_\tau^k \right\} + \right. \\
\left. \left[ (\dd_{n\Omega}^k + \dd_{nz}^k) f_x \delta\uu_s^k)^{\rm T} (\cc_{np}^k \dd_p^k F_\tau \uu_\tau^k  + \cc_{nn}^k (\dd_{n\Omega}^k + \dd_{nz}^k) F_\tau \uu_\tau^k) \right] \right\}~\mathrm{d}\Omega_k ~\mathrm{d}z = \\ \sum\limits_{k=1}^{N_k} \int\limits_{\Omega_k} \int\limits_{A_k} \rho^k \delta \uu_s^{k^{\rm T}} \ddot{\uu}^k~\mathrm{d}\Omega_k~\mathrm{d}z + \sum\limits_{k=1}^{N_k} \delta \mathbf{L}_e^k
\end{split}
\end{equation}
After integration by parts, the governing equations for the plate are obtained:
\begin{equation}
\kk_{uu}^{k\tau s} \uu_{\tau}^k = \mathbf{P}_{u \tau}^k
\end{equation}
and in the case of free vibrations, we have:
\begin{equation}
\kk_{uu}^{k\tau s} \uu_{\tau}^k = \mathbf{M}^{k \tau s} \ddot{\uu}_\tau^k
\end{equation}
where the fundamental nucleus $\kk_{uu}^{k \tau s}$ is:
\begin{equation}
\kk_{uu}^{k \tau s} = \left[ (-\dd_p^k)^{\rm T} ( \cc_{pp}^k \dd_p^k + \cc_{pn}^k (\dd_{n\Omega}^k + \dd_{nz}) + (-\dd_{n\Omega}^k + \dd_{nz}^k)^{\rm T} (\cc_{np}^k \dd_p^k + \cc_{nn}^k (\dd_{n\Omega}^k + \dd_{nz}^k)) \right] F_\tau F_s
\label{eqn:stifffundanuclei}
\end{equation}
and $\mathbf{M}^{k \tau s}$ is the fundamental nucleus for the inertial term given by:
\begin{equation}
M_{ij}^{k \tau s} = \left\{ \begin{array}{cc} \rho^k F_\tau F_s & \textup{if} \hspace{1cm} i = j \\ 0 & \textup{if} \hspace{1cm} i \neq j \end{array} \right.
\label{eqn:massfundanuclei}
\end{equation}
where $\mathbf{P}_{u \tau}^k$ are variationally consistent loads with applied pressure. For more detailed derivation and for the explicit form of the fundamental nuclei, interested readers are referred to~\cite{carrerademasi2002,carrerademasi2002a}. The work done by the applied non-conservative loads is:
\begin{equation}
W(\boldsymbol{\delta}) = \int_{\Omega} \Delta p w ~\rmd \Omega
\label{eqn:aerowork}
\end{equation}
where $\Delta p$ is the aerodynamic pressure. The static aerodynamic pressure based on first-order, high Mach number approximation to linear potential flow is given by:
\begin{equation}
\Delta p = \frac{\rho_a U_a^2}{\sqrt{M_\infty^2 - 1}} \left[ \frac{\partial w}{\partial x} \cos \theta^\prime + \frac{\partial w}{\partial y} \sin \theta^\prime \right]
\label{eqn:aeropressure}
\end{equation}
where $\rho_a, U_a, M_\infty$ and $\theta^\prime$ are the free stream air density, velocity of air, Mach number and flow angle, respectively. 
\begin{equation}
\left[ \left( \KK + \lambda \overline{\mathbf{A}}\right) - \overline{\kappa} \bfm\right] \boldsymbol{\delta} = \mathbf{0}
\label{eqn:finaldiscre}
\end{equation}
where $\KK$ is the stiffness matrix, $\overline{\mathbf{A}}$ is the aerodynamic matrix and $\bfm$ is the mass matrix. The eigenvalue $\overline{\kappa} = -\omega^2 - g_\tau \omega/(\rho h)$ includes the contribution of aerodynamic damping. \Eref{eqn:finaldiscre} is solved for eigenvalues for a given value of $\lambda$. In the absence of aerodynamic damping, i.e., when $\lambda = $ 0, the eigenvalue, $\omega$ is real and positive, since the stiffness matrix and mass matrix are symmetric and positive definite. However, the aerodynamic matrix $\overline{\mathbf{A}}$ is unsymmetric and hence complex eigenvalues $\omega$ are expected for $\lambda >$ 0. As $\lambda$ increases monotonically from zero, two of these eigenvalues will approach each other and become complex conjugates. In this study, $\lambda_{cr}$ is considered to be the value of $\lambda$ at which the first coalescence occurs. In the presence of aerodynamic damping, the eigenvalues $\overline{\kappa}$, in \Eref{eqn:finaldiscre} becomes complex with increasing value of $\lambda$. The corresponding frequency can be written as:
\begin{equation}
\overline{\kappa} = -\omega^2 - g_\tau \omega/(\rho h) = \overline{\kappa}_R - i \overline{\kappa}_I
\end{equation}
where the subscripts $R$ and $I$ refer to the real and the imaginary part of the eigenvalue. The flutter boundary is reached $(\lambda = \lambda_{cr})$, when the frequency $\omega$ becomes pure imaginary number, i.e., $\omega = i \sqrt{\overline{\kappa}_R}$ at $g_\tau = \overline{\kappa}_I/\sqrt{\overline{\kappa}_R}$. In practice, the value of $\lambda_{cr}$ is determined from a plot of $\omega_R$ vs $\lambda$ corresponding to $\omega_R = $ 0. 

\section{Element description} \label{edesc}
The plate element employed in this study is a $\mathcal{C}^0$ continuous element and according to the isoparametric description, the components of each displacement unknown $\uu_\tau$ are expressed as:
\begin{equation}
\uu_\tau = N_I \qq_{\tau I}, \hspace{0.5cm} I = 1,2,\cdots,N_n
\label{eqn:unifieddisp}
\end{equation}
where $N_I$ are the standard finite element shape functions. By introducing the unified formulation for the displacements, given by \Eref{eqn:unifieddisp} into the strain-displacement relations (see \Eref{eqn:strainDisp}), we have:
\begin{align}
\bveps_{pG}^k &= \dd_p^k (F_\tau \uu_\tau^k) = \dd_p^k (F_{\tau} N_I) \qq_{\tau I}^k \nonumber \\
\bveps_{nG}^k &= (\dd_{n\Omega}^k + \dd_{nz}^k) (F_\tau \uu_\tau^k) = \dd_{n\Omega}^k (F_\tau N_I) \qq_{\tau I}^k + F_{\tau,z} N_I \qq_{\tau I}^k
\label{eqn:straindiscretize}
\end{align}
Upon substituting \Erefs{eqn:unifieddisp} and (\ref{eqn:straindiscretize}) into \Erefs{eqn:stifffundanuclei} and (\ref{eqn:massfundanuclei}), we can compute the stiffness matrix $\KK$, the aerodynamic matrix $\overline{\mathbf{A}}$ and the mass matrix $\bfm$ of the system. The formulation is implemented in MATLAB\textsuperscript{\textregistered} and the solution to the \Eref{eqn:finaldiscre} is computed from a standard eigenvalue algorithm.

\paragraph{Shear locking} If the interpolation functions given for a QUAD-4 are used directly to interpolate the unknown displacement fields in deriving the shear strains $(\gamma_{xz}, \gamma_{yz})$ and the membrane strains $(\bveps_{pG})$, the element will lock and show oscillations in the shear and the membrane stresses. The oscillations are due to the fact that the derivative functions of the out-of plane displacement do not match that of the rotations in the shear strain definition. To alleviate the locking phenomenon, the terms corresponding to the derivative of the out-of plane displacement must be consistent with the rotation terms. In this study, field redistributed shape functions are used to alleviate shear locking. The field consistency requires that the transverse shear strains and the membrane strains must be interpolated in a consistent manner. If the element has edges which are aligned with the coordinate system $(x,y)$, the terms in shear strains $(\gamma_{xz}, \gamma_{yz})$ are approximated by the following substitute shape functions~\cite{somashekarprathap1987}:
\begin{eqnarray}
\tilde{N}_{1}(\eta) = \frac{1}{4} \left[ \begin{array}{cccc} 1-\eta & 1-\eta & 1+\eta &  1+\eta \end{array} \right] \nonumber \\
\tilde{N}_{2}(\xi) = \frac{1}{4} \left[ \begin{array}{cccc} 1-\xi  & 1 +\xi & 1+\xi  & 1-\xi \end{array} \right].
\label{eqn:fieldredistribute}
\end{eqnarray}
Note that, no special integration rule is required for evaluating the shear terms. A numerical integration based on the 2 $\times$ 2 Gaussian rule is used to evaluate all the terms.

\section{Numerical Results}
\label{numres}
In this section, we present the critical aerodynamic pressure and the critical frequency of laminated composite plates immersed in a supersonic flow using 4-noded quadrilateral element and unified formulation. In this study, we use a hybrid displacement assumption, where the in-plane displacements $u$ and $v$ are expressed as sinusoidal expansion in the thickness direction, and the transverse displacement, $w$ is quadratic in the thickness direction. We refer to this theory as SINUS-W2. The displacements are expressed as:
\begin{align}
u(x,y,z,t) &= u_o(x,y,t) + zu_1(x,y,t) + \sin \left( \frac{\pi z}{h} \right) u_2(x,y,t) \nonumber \\
v(x,y,z,t) &= v_o(x,y,t) + zv_1(x,y,t) + \sin \left( \frac{\pi z}{h}\right) v_2(x,y,t) \nonumber \\
w(x,y,z,t) &= w_o(x,y,t) + zw_1(x,y,t) + z^2 w_2(x,y,t)
\end{align}
where $u_o, v_o$ and $w_o$ are translations of a point at the middle-surface of the plate, $w_2$ is higher order translation, and $u_1, v_1, u_3$ and $v_3$ denote rotations~\cite{touratier1991} and considers a quadratic variation of the transverse displacement $w$ allowing for through-the-thickness deformations. Both simply supported and clamped boundary conditions are considered in this study and the influence of the flow direction is also studied. In all cases, we present the non dimensionalized critical aerodynamic pressure, $\lambda_{cr}$ and critical frequency $\omega_{cr}$ as, unless specified otherwise:
\begin{eqnarray}
\omega^\ast_{cr} = \omega_{cr} a^2 \sqrt{ \frac{\rho h}{D}} \nonumber \\
\lambda^\ast_{cr} = \lambda_{cr} \frac{a^3}{D}
\label{eqn:nondimfreq}
\end{eqnarray}
where $D = {E_2 h^3 \over 12(1-\nu^2)}$ is the bending rigidity of the plate, $E_2, \nu$ are the Young's modulus and Poisson's ratio and $\rho$ is the mass density. 

Before proceeding with the detailed numerical study, the formulation developed herein is validated against available results pertaining to the critical aerodynamic pressure and the critical frequency for a 5-layered laminated square plate. Table \ref{tab:meshconvefreq} presents the convergence of the first three fundamental frequencies with mesh size. A structured mesh of 30$\times$30 is found to be adequate for this study. It can be seen that the results from the present formulation are in good agreement with those in the literature. Further numerical studies are performed with a structured quadrilateral mesh. Table \ref{tab:meshconveflutter} presents the convergence of the flutter bounds for a clamped square laminated plate with $a/h=$ 100. The flutter bounds are computed for both the cases: without damping and with damping. It is seen that the results agree well with those in the literature.

\begin{table}[htbp]
\centering
\renewcommand{\arraystretch}{1.2}
\caption{Convergence of the non-dimensional natural frequencies $\Omega = \omega \frac{a^2}{\pi^2 h} \sqrt{ \frac{\rho}{E_2}} $ of a 5-layered laminated square plate $[45^\circ/-45^\circ/45^\circ/-45^\circ/45^\circ]$ with $E_L/E_T=$ 40, $G_{LT}/E_{T}=$ 0.6, $G_{TT}/E_T=$ 0.5, $\nu_{LT}=$0.25. }
\begin{tabular}{llrrr}
\hline
Mesh && \multicolumn{3}{c}{Modes} \\
\cline{3-5}
&& Mode 1 & Mode 2 & Mode 3 \\
\hline
5$\times$5	&& 2.5593	& 5.7514 & 7.2116 \\
10$\times$10	&& 2.4571 &	5.1479 & 6.3702 \\
14$\times$14	&& 2.4413	& 5.0508 & 6.2475 \\
20$\times$20	&& 2.4316	& 5.0008 & 6.1809 \\
30$\times$30	&& 2.4254	& 4.9746 & 6.1435 \\
Ref.~\cite{singhaganapathi2005}	 && 2.4343 & 4.9854 & 6.1823 \\
Ref.~\cite{wang1997} && 2.4339 & 4.9865 & 6.1818 \\
\hline
\end{tabular}
\label{tab:meshconvefreq}
\end{table}

\begin{table}[htbp]
\centering
\renewcommand{\arraystretch}{1.2}
\caption{Convergence of natural frequencies and flutter bounds for $[(0^\circ/90^\circ)_{2\rm{s}}]$ boron/epoxy clamped laminate with $a/b=$ 1.0 and $a/h=$ 100. }
\begin{tabular}{lrrcrrrrr}
\hline
Mesh & \multicolumn{2}{c}{Modes} && \multicolumn{2}{c}{Without damping} && \multicolumn{2}{c}{With Damping}\\
\cline{2-3}\cline{5-6}\cline{8-9}
& Mode 1 & Mode 2 && $\lambda^\ast_{cr}$ & $\omega^\ast_{cr}$ && $\lambda^\ast_{cr}$ & $\omega^\ast_{cr}$\\
\hline
5$\times$5	& 26.5697 & 59.0459  && 866.60 & 58.78 && 883.44 & 59.23  \\
10$\times$10	& 24.1446 & 45.5384  && 513.48 & 48.37 && 530.31 & 49.02\\
14$\times$14	& 23.7904 & 43.9233  && 479.88 & 47.24 && 496.72 & 47.93 \\
20$\times$20	& 23.6064 & 43.1125  && 464.26 & 46.71 && 481.09 & 47.42 \\
30$\times$30 & 23.5094 & 42.6922 && 455.66 & 46.41 && 472.50 & 47.14 \\
Integral equation method~\cite{srinivasanbabu1987} & 23.33 &  53.77 && - & - && 446.36	 & 46.09 \\
Series solution.~\cite{srinivasanbabu1987} & 23.63 &  53.76 && - & - && 474.60 & 47.19 \\
Classical lamination theory~\cite{leecho1991} & 23.34 & 42.30 && - & - && 471.16 & 46.89 \\
\hline
\end{tabular}
\label{tab:meshconveflutter}
\end{table}

\begin{figure}
\centering
\setlength\figureheight{8cm} 
\setlength\figurewidth{10cm}
%
%
%
%
\begin{tikzpicture}

\begin{axis}[%
width=\figurewidth,
height=\figureheight,
scale only axis,
xmin=0,
xmax=180,
xlabel={$\text{Flow angle, }\theta^\prime$},
ymin=120,
ymax=300,
ylabel={Critical Pressure, $\lambda^\ast_{\rm {cr}}$},
legend style={draw=black,fill=white,legend cell align=left}
]
\addplot [
color=black,
solid,
mark=square,
mark options={solid}
]
table[row sep=crcr]{
0 258.0078\\
10 260.3516\\
20 268.4766\\
30 262.0703\\
40 220.3516\\
50 188.0078\\
60 164.2578\\
70 153.3203\\
80 147.0703\\
90 144.7266\\
100 147.0703\\
110 153.3203\\
120 164.2578\\
130 188.0078\\
140 220.3516\\
150 262.0703\\
160 268.4766\\
170 260.3516\\
180 258.0078\\
};
\addlegendentry{$\lambda^\ast_{\rm {cr}}$};

\end{axis}

\begin{axis}[%
width=\figurewidth,
height=\figureheight,
scale only axis,
xmin=0,
xmax=180,
ymin=20,
ymax=40,
hide x axis,
axis y line*=right,
ylabel={Critical Frequency, $\omega^\ast_{\rm {cr}}$},
ylabel near ticks,
legend style={draw=black,fill=white,legend cell align=left},
legend pos = south west
]

\addplot [
color=black,
solid,
mark=o,
mark options={solid}
]
table[row sep=crcr]{
0 29.1219\\
10 29.4612\\
20 30.5337\\
30 30.2253\\
40 27.0337\\
50 25.1255\\
60 23.9001\\
70 23.4254\\
80 23.168\\
90 23.0678\\
100 23.168\\
110 23.4254\\
120 23.9001\\
130 25.1255\\
140 27.0337\\
150 30.2253\\
160 30.5337\\
170 29.4612\\
180 29.1219\\
};
\addlegendentry{$\omega^\ast_{\rm {cr}}$};

\end{axis}
\end{tikzpicture}%
\caption{Effect of flow angle $\theta^\prime$ on the critical aerodynamic pressure and critical frequency for a simply supported square laminated plate with $[(0^\circ/90^\circ)_{\rm{2s}})]$ and $a/b=$ 1, $a/h=$ 100.}
\label{fig:thratioonlambdaomega}
\end{figure}
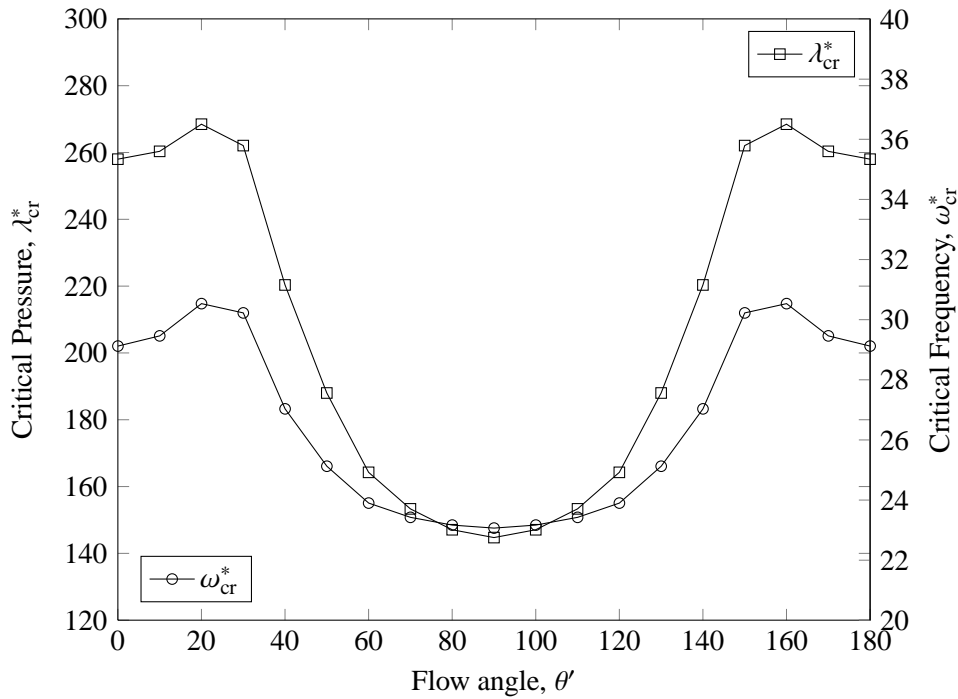

Next, the influence of the boundary conditions, the plate thickness and the direction of the flow on the flutter bounds are numerically investigated. Table \ref{tab:ahbcinfluence} presents the influence of the boundary conditions and the plate thickness for a square laminated plate with the following stack up sequence $[(0^\circ/90^\circ)_{\rm{2s}})]$ immersed in a normal flow. It is seen that with increasing plate thickness, the flutter bounds decreases as expected. Also, the effect of damping is to increase the critical aerodynamic pressure and the critical frequency. The effect of the boundary conditions on the flutter bounds is also seen in Table \ref{tab:ahbcinfluence}. The influence of the plate aspect ratio $a/b$ and the flow angle $\theta^\prime$ is shown in \frefs{fig:abratioonlambdaomega} - \ref{fig:thratioonlambdaomega}. The flow is considered to be along the x-direction. It is seen from \fref{fig:abratioonlambdaomega} that as the plate aspect ratio $(a/b)$ increases, the flutter bounds decreases for the stack up sequence considered here. The effect of the flow angle $\theta^\prime$ on the flutter bounds of a square plate is shown in \fref{fig:thratioonlambdaomega}. Increasing the flow angle, the flutter bounds, viz., the critical aerodynamic pressure $\lambda^\ast_{\rm cr}$ and the critical frequency $\omega^\ast_{\rm cr}$, initially increases until $\theta^\prime=$ 20$^\circ$. Upon further increasing in the flow angle, the flutter bounds decreases monotonically until it reaches a minimum at $\theta^\prime=$ 90$^\circ$. The flutter behaviour is symmetric with respect to the flow angle $\theta^\prime=$ 90$^\circ$.

\begin{table}[htbp]
\centering
\renewcommand{\arraystretch}{1.2}
\caption{Influnence of the plate thickness $a/h$ and the support conditions on the flutter characteristics of laminated plates immersed in a supersonic flow with $a/b=$ 1 and flow angle $\theta^\prime=$ 0$^\circ$.}
\begin{tabular}{lrrcrrrrrr}
\hline
Boundary & $a/h$& \multicolumn{2}{c}{Modes} && \multicolumn{2}{c}{Without damping} && \multicolumn{2}{c}{With Damping}\\
\cline{3-4}\cline{6-7}\cline{9-10}
& & Mode 1 & Mode 2 && $\lambda^\ast_{cr}$ & $\omega^\ast_{cr}$ && $\lambda^\ast_{cr}$ & $\omega^\ast_{cr}$\\
\hline
\multirow{2}{*}{CCCC} & 100 & 23.5094 & 42.6922 && 455.66 & 46.41 && 472.50 & 47.14 \\
& 10 & 15.7281 & 25.2179 && 139.26 & 27.05 && 145.51 & 27.74 \\
\cline{2-10}
\multirow{2}{*}{SSSS} & 100 & 10.9018 & 26.4822 && 251.76 & 28.67 && 258.01 & 29.12 \\
& 10 & 9.6205 & 16.5652 && 154.88 & 24.39 && 166.14 & 25.62 \\
\hline
\end{tabular}
\label{tab:ahbcinfluence}
\end{table}

\begin{figure}
\centering
\newlength\figureheight 
\newlength\figurewidth 
\setlength\figureheight{8cm} 
\setlength\figurewidth{10cm}
%
%
%
%
\begin{tikzpicture}

\begin{axis}[%
width=\figurewidth,
height=\figureheight,
scale only axis,
xmin=0.4,
xmax=2,
xlabel={a/b},
ymin=240,
ymax=400,
ylabel={Critical Pressure, $\lambda^\ast_{\rm {cr}}$},
legend style={draw=black,fill=white,legend cell align=left}
]
\addplot [
color=black,
solid,
mark=square,
mark options={solid}
]
table[row sep=crcr]{
0.4 350.1953\\
0.6 285.6641\\
0.8 270.0391\\
1 258.0078\\
1.2 254.1016\\
1.4 251.7578\\
1.6 250.1953\\
1.8 249.4141\\
2 248.6328\\
};
\addlegendentry{$\lambda^\ast_{\rm {cr}}$};

\end{axis}

\begin{axis}[%
width=\figurewidth,
height=\figureheight,
scale only axis,
xmin=0.4,
xmax=2,
ymin=10,
ymax=52,
hide x axis,
axis y line*=right,
ylabel={Critical Frequency, $\omega^\ast_{\rm {cr}}$},
ylabel near ticks,
legend style={draw=black,fill=white,legend cell align=left},
legend pos = south west
]

\addplot [
color=black,
solid,
mark=o,
mark options={solid}
]
table[row sep=crcr]{
0.4 50.2302\\
0.6 34.5076\\
0.8 30.6863\\
1 29.1219\\
1.2 28.5282\\
1.4 28.2211\\
1.6 28.0066\\
1.8 27.9416\\
2 27.8643\\
};
\addlegendentry{$\omega^\ast_{\rm {cr}}$};

\end{axis}
\end{tikzpicture}%
\caption{Influence of the plate aspect ratio $a/b$ on the flutter parameters, viz., the critical pressure $\lambda^\ast_{\rm cr}$ and the critical frequency $\omega^\ast_{cr}$ for a simply supported square plate with $a/h=$ 100. The flow is normal to the plate, i.e., flow angle $\theta^\prime=$ 0$^\circ$.}
\label{fig:abratioonlambdaomega}
\end{figure}
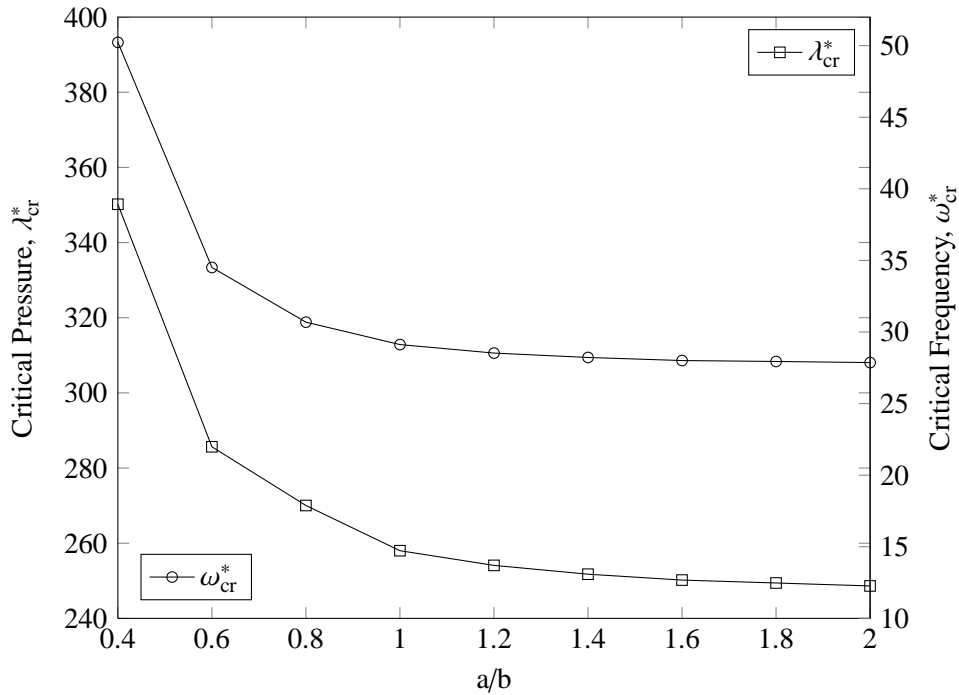

\section{Conclusions}
The flutter characteristics of laminated composites immersed in a supersonic flow has been analyzed within the framework of unified formulation. The plate kinematics is based on sinusoidal theory and a quadratic variation of the transverse displacement through the thickness is considered. A shear flexible four noded quadrilateral plate element was used to discretize the domain and the aerodynamic force is accounted for assuming the first-order Mach number approximation potential flow theory. The results from the present formulation are in very good agreement with the results available in the literature. The influence of the plate aspect ratio and the flow angle are numerically studied.

\section*{Acknowledgements} 
S Natarajan would like to acknowledge the financial support of the School of Civil and Environmental Engineering, The University of New South Wales for his research fellowship for the period September 2012 onwards. 

\section*{References}
\bibliographystyle{elsarticle-num}
\bibliography{Flut}

\end{document}